\documentclass[12pt,twoside,a4paper]{article}
\usepackage[russian]{babel}
\usepackage{amsmath,amsfonts,amscd,amssymb,latexsym}
\usepackage[dvips]{graphicx}

\begin{document}

\sloppy
\begin{center} 
{\large\bf GENERALIZED KANTOR DOUBLE}\\

\hspace*{8mm}

{\large\bf Ivan Kaygorodov}

\

{\it 
Sobolev Inst. of Mathematics\\ 
Novosibirsk, Russia\\
kib@math.nsc.ru\\}

\

\

\end{center}

\begin{center} {\bf Abstract: }\end{center}                                                                    
We find necessary and sufficient conditions for a generalized Kantor double to be Jordan. 
We also describe $\delta$-superderivations of a generalized Kantor double whose even part is prime.

\

{\bf Key words:} Kantor double, Jordan superalgebra, $\delta$-superderivation. 

\

\section{Предварительные сведения}

Пусть $F$ --- поле характеристики $p \neq 2$. Супералгеброй над полем $F$ называется алгебра $A$, 
такая что $A=A_0\oplus A_1$ с условием $A_i A_j \subseteq  A_{i+j(mod2)}.$
Элементы супералгебры $A=A_0 \oplus A_1$ из множества $A^{*}=A_0 \cup A_1$ будем называть однородными.
Для однородного элемента $x$ супералгебры $A$ будем считать $p(x)=i,$ если $x\in A_i.$

Алгебра $A$ над полем $F$
называется йордановой, если она удовлетворяет тождествам
$$xy=yx, (x^{2}y)x=x^{2}(yx).$$

Пусть  $G$ --- алгебра Грассмана над $F$, заданная
образующими $1,\xi_{1},\ldots ,\xi_{n},\ldots $ и определяющими
соотношениями: $\xi_{i}^{2}=0, \xi_{i}\xi_{j}=-\xi_{j}\xi_{i}.$ Элементы
$1, \xi_{i_{1}}\ldots \xi_{i_{k}} (i_{1}< \ldots
<i_{k})$ образуют базис алгебры $G$ над $F$. Обозначим через
$G_{0}$ и $G_{1}$ подпространства, порожденные,
соответственно, произведениями четной и нечетной длины; тогда
$G$ представляется в виде прямой суммы этих подпространств:
$G = G_{0}\oplus G_{1}$, при этом справедливы
соотношения $G_{i}G_{j} \subseteq G_{i+j(mod 2)},
i,j=0,1.$ 

Под суперпространством мы понимаем $\mathbb{Z}_{2}$-градуированное
пространство. На пространстве $End(A)$ эндоморфизмов супералгебры $A=A_0+A_1$ зададим структуру супералгебры, 
таким образом, что четными элементами будем считать те эндоморфизмы, которые инвариатны на $A_0$ и $A_1$,
а нечетными элементами будем считать такие эндоморфизмы $\phi$, что $\phi(A_i) \subseteq A_{i+1}.$

Для супералгебры $A$ подалгебра 
$$G(A)=G_{0} \otimes A_{0} + G_{1} \otimes A_{1}$$
в тензорном произведении $G \otimes A$ называется грассмановой оболочкой супералгебры $A$.

Если $\Omega$ --- некоторое многообразие алгебр над $F$.
Супералгебра $A$ называется $\Omega$-супералгеброй,
если $G(A) \in \Omega$. Таким образом, супералгебра $A$ является йордановой супералгеброй, 
если ее грассманова оболочка $G(A)$ является йордановой алгеброй.

{\bf Дубль Кантора \cite{kant}.}
Пусть $\Gamma=\Gamma_0 \oplus \Gamma_1$ --- ассоциативная суперкоммутативная супералгебра с единицей 1 и $\{, \}: 
\Gamma \times \Gamma \rightarrow \Gamma$ --- суперкососимметрическое билинейное отображение, которое мы будем называть скобкой. 
По супералгебре $\Gamma$ и скобке $\{, \}$ можно построить супералгебру $J(\Gamma,\{,\})$. 
Рассмотрим $J(\Gamma, \{,\})=\Gamma \oplus \Gamma x$ --- прямую сумму пространств, 
где $\Gamma x$ --- изоморфная копия пространства $\Gamma.$ 
Считаем, что $D(a)=\{a,1\}.$ Пусть $a,b$ --- однородные элементы из $\Gamma$. 
Тогда операция умножения $\cdot$ на $J(\Gamma, \{, \})$ определяется формулами
$$a\cdot b=ab, a\cdot bx=(ab)x, ax\cdot b=(-1)^{p(b)}(ab)x, ax \cdot bx = (-1)^{p(b)}\{a, b\}.$$

Положим $A=\Gamma_0+\Gamma_1 x, M = \Gamma_1+\Gamma_0 x.$ Тогда $J(\Gamma, \{, \})=A\oplus M$ --- $\mathbb{Z}_2$-градуированная алгебра.

Для унитальной супералгебры скобка $\{,\}$ называется йордановой, если при однородных элементах $f_i,g_i,h_i \in \Gamma_i$ выполняются 
следующие соотношения
\begin{eqnarray}\label{jo1}\{f_i,g_jh_k\}=\{f_i,g_j\}h_k+(-1)^{ij}g_j\{f_i,h_k\}-D(f_i)g_jh_k, \end{eqnarray}
$$\{f_i,\{g_j,h_k\}\}=\{\{f_i,g_j\},h_k\}+(-1)^{ij}\{g_j,\{f_i,h_k\}\}+$$
\begin{eqnarray}\label{jo2}D(f_i)\{g_j,h_k\}+(-1)^{ji}D(g_j)\{h_k,f_i\}+(-1)^{k(j+i)}D(h_k)\{f_i,g_j\}. \end{eqnarray}
В дальнейшем, элементы $k_i,f_i,g_i,h_i$ мы будем всегда считать однородными.

Хорошо известно \cite{KingMac2, KingMac}, что  супералгебра $J(\Gamma, \{, \})$ йорданова тогда и только тогда, 
когда скобка $\{,\}$ является йордановой. В силу йордановости супералгебры $J(\Gamma, \{, \})$ получаем, 
что $D: a \rightarrow \{a,1\}$ --- дифференцирование супералгебры $\Gamma$.

Если $D$ --- нулевое дифференцирование, то $\{,\}$ является скобкой Пуассона, т.е. 
$$\{a,bc\}=\{a,b\}c+(-1)^{p(a)p(b)}b\{a,c\}$$ и $\Gamma$ --- супералгебра Ли относительно операции $\{,\}$.  Произвольная скобка Пуассона является йордановой скобкой \cite{kant2}.

\medskip

Хорошо известно \cite{KingMac2, KingMac}, что йорданова супералгебра
$J=\Gamma+\Gamma x,$ полученная с помощью процесса удвоения
Кантора, будет являться простой тогда и только тогда, когда
$\Gamma$ не имеет ненулевых идеалов $B$ с условием $\{\Gamma,B\}
\subseteq B.$ 

\medskip

\textbf{Супералгебра векторного типа $J(\Gamma, D)$.} Пусть $\Gamma$ --- ассоциативная суперкоммутативная супералгебра 
с ненулевым четным дифференцированием $D$. Определим на $\Gamma$ скобку $\{,\}$ полагая $\{a,b\}=D(a)b-aD(b).$ 
Тогда скобка $\{,\}$ --- йорданова скобка. Полученную супералгебру $J(\Gamma, \{,\})$ будем обозначать как $J(\Gamma,D)$.

\medskip

\section{Йордановость обобщенного дубля Кантора}

В данной части работы мы рассмотрим обобщенный дубль Кантора. Мы ослабим условие унитальности и суперкоммутативности 
супералгебры $\Gamma$ и построим супералгебру $J(\Gamma, \{,\})$ по аналогии с выше приведенной конструкцией. Если супералгебра 
$\Gamma$ обладает дифференцированием $D$, то, задав скобку $\{,\}$ по правилу $\{a,b\}=D(a)-aD(b)$, мы получим супералгебру 
векторного типа (не обязательно унитальную).

Легко заметить, что условие простоты супералгебры $J(\Gamma, \{,\})$ влечет отсутствие идеалов $I$ в супералгебре $\Gamma$ c 
условием $\{I,\Gamma\} \subseteq I.$ В противном случае, в супералгебре $J(\Gamma, \{,\})$ мы имели бы ненулевой градуированный
идеал $I+Ix$.

Скобка $\{, \}$, определенная на супералгебре $\Gamma,$  называется йордановой, если 

$$(-1)^{(i+j)l}\{\{f_{i},h_k\}g_j,k_l\}+(-1)^{(k+j)i}\{\{h_k,k_l\}g_j,f_i\}+$$
$$(-1)^{(l+j)k}\{\{k_l,f_i\}g_j,h_k\}=(-1)^{(i+j)l}\{f_i,h_k\}\{g_j,k_l\}+$$
\begin{eqnarray}\label{jorskob1}(-1)^{(k+j)i}\{h_k,k_l\}\{g_j,f_i\}+(-1)^{(l+j)k}\{k_l,f_i\}\{g_j,h_k\},\end{eqnarray}

$$(-1)^{(k+j)i}(\{h_kk_l,g_j\}f_i-h_kk_l\{g_j,f_i\})=$$
\begin{eqnarray}\label{jorskob2}(-1)^{(l+j)k}(\{k_lf_i,g_j\}h_k-k_lf_i\{g_j,h_k\}),\end{eqnarray}

$$(-1)^{(i+j)l}(\{f_ih_kg_j,k_l\}-f_ih_k\{g_j,k_l\})=$$
$$(-1)^{(k+j)i}(\{h_kk_l,g_j\}f_i-\{h_kk_l,g_jf_i\})+$$
\begin{eqnarray}\label{jorskob3}
(-1)^{(l+j)k}(\{k_lf_i,g_j\}h_k-\{k_lf_i,g_jh_k\}).\end{eqnarray}

{\bf Теорема 2.1.} Обобщенный дубль Кантора $J(\Gamma, \{,\})$ является йордановой супералгеброй тогда и только тогда, когда 
скобка $\{, \}$ является йордановой и супералгебра $\Gamma$ --- суперкоммутативна.

\medskip 

{\bf Доказательство.} Известно, что йорданова супералгебра удовлетворяет тождествам 
\begin{eqnarray}\label{jor1}{\bf f_i \cdot g_j} - (-1)^{ij} {\bf g_j \cdot f_i}=0,\end{eqnarray}
$(-1)^{(i+j)l}{\bf [f_i \cdot h_k, g_j,k_l]}+$
\begin{eqnarray}\label{jor2}(-1)^{(k+j)i}{\bf [h_k \cdot k_l, g_j, f_i]} +(-1)^{(l+j)k}{\bf [k_l\cdot f_i, g_j, h_k]}=0,\end{eqnarray}
где $[f,g,h]=(f \cdot g)\cdot h - f \cdot (g \cdot h)$ --- обычный неградуированный ассоциатор и ${\bf f_i}=f_i+ f_{i+1}x$ 
--- базисные элементы $J_i=F_i+F_{i+1}x.$ 

Легко заметить, что супералгебра $J(\Gamma, \{,\})$ удовлетворяет условию (\ref{jor1}) тогда и только тогда, когда $\Gamma$ 
--- суперкоммутативная супералгебра. 

Для установления эквивалентности между тождеством (\ref{jor2}) и соотношениями йордановой скобки (\ref{jorskob1}-\ref{jorskob3})
 нам достаточно проверить эквивалентность 
на однородных элементах супералгебры $J(\Gamma, \{,\}).$ В дальнейшем, считаем, что $f_t,g_t,h_t,k_t \in \Gamma_t.$ Поскольку $\Gamma$ --- ассоциативна, то 
$$[f_i \cdot h_k, g_j,k_l]=0,[f_ix \cdot h_k,g_j,k_l]=0, [f_i \cdot h_kx,g_j,k_l]=0, $$
$$[f_i \cdot h_k,g_jx,k_l]=0,[f_i \cdot h_k,g_j,k_lx]=0.$$
Таким образом, достаточно рассмотреть случаи, когда в тождестве (\ref{jor2}) среди элементов ${\bf f_i,g_j, h_k, k_l}$ 
два или более являются элементами $\Gamma x.$ Проведем последовательные вычисления.

\medskip

$0=
(-1)^{(i+j)l}[f_{i+1}x \cdot h_{k+1}x, g_{j+1}x, k_{l+1}x]+\\
(-1)^{(k+j)i}[h_{k+1}x \cdot k_{l+1}x, g_{j+1}x, f_{i+1}x]+ 
(-1)^{(l+j)k}[k_{l+1}x\cdot f_{i+1}x, g_{j+1}x, h_{k+1}x]=\\
(-1)^{(i+j)l+k+l}(\{\{f_{i+1},h_{k+1}\}g_{j+1},k_{l+1}\}-\{f_{i+1},h_{k+1}\}\{g_{j+1},k_{l+1}\})+\\
(-1)^{(k+j)i+l+i}(\{\{h_{k+1},k_{l+1}\}g_{j+1},f_{i+1}\}-\{h_{k+1},k_{l+1}\}\{g_{j+1},f_{i+1}\})+\\
(-1)^{(l+j)k+i+k}(\{\{k_{l+1},f_{i+1}\}g_{j+1},h_{k+1}\}-\{k_{l+1},f_{i+1}\}\{g_{j+1},h_{k+1}\}).$

Легко заметить, что мы получаем аналог тождества (\ref{jorskob1}).

\medskip

$0=
(-1)^{(i+j)l}[f_{i+1}x \cdot h_{k+1}x, g_{j+1}x,k_l]+$
\begin{eqnarray}\label{anal2.4}(-1)^{(k+j)i}[h_{k+1}x\cdot k_l, g_{j+1}x, f_{i+1}x]+
(-1)^{(l+j)k}[k_l\cdot f_{i+1}x, g_{j+1}x,h_{k+1}x]=\end{eqnarray}
$(-1)^{(i+j)l+k+l+1}(\{f_{i+1},h_{k+1}\}g_{j+1}k_l-\{f_{i+1},h_{k+1}\}g_{j+1}k_l)x+\\
(-1)^{(k+j)i+l+j+1}(\{h_{k+1}k_l,g_{j+1}\}f_{i+1}-h_{k+1}k_l\{g_{j+1},f_{i+1}\})x+\\
(-1)^{(l+j)k+j+1}(\{k_lf_{i+1},g_{j+1}\}h_{k+1}-k_lf_{i+1}\{g_{j+1},h_{k+1}\})x.$

Заметим, что мы имеем аналог тождества (\ref{jorskob2}).

\medskip 

Вычисляя, видим что $[fx \cdot hx,g, kx]=0,$ следовательно, при 
подстановке ${\bf f_i}=f_{i+1}x, {\bf g_j}=g_{j}, {\bf h_k}=h_{k+1}x, {\bf k_l}=k_{l+1}x$ в 
соотношение (\ref{jor2}), в правой части мы имеем нулевое выражение.

\medskip

$0=
(-1)^{(i+j)l}[f_{i+1}x \cdot h_k, g_{j+1}x, k_{l+1}x]+\\
(-1)^{(k+j)i}[h_k \cdot k_{l+1}x, g_{j+1}x, f_{i+1}x]+
(-1)^{(l+j)k}[k_{l+1}x \cdot f_{i+1}x, g_{j+1}x, h_k]$

Отметим, что при замене $f_i \mapsto h_k, h_k \mapsto k_l, k_l \mapsto f_i$ мы получим выражение (\ref{anal2.4}), которое эквивалентно
(\ref{jorskob2}).

\medskip 

$0=
(-1)^{(i+j)l}[f_i \cdot h_{k+1}x, g_{j+1}x, k_{l+1}x]+\\
(-1)^{(k+j)i}[h_{k+1}x\cdot k_{l+1}x, g_{j+1}x, f_i]+
(-1)^{(l+j)k}[k_{l+1}x \cdot f_i, g_{j+1}x, h_{k+1}x]$

Элементарно заметить, что при замене $f_i \mapsto k_l, k_l \mapsto h_k, h_k \mapsto f_i$ мы получим выражение (\ref{anal2.4}), которое эквивалентно
(\ref{jorskob2}).

\medskip 

$0=
(-1)^{(i+j)l}[f_i\cdot h_k, g_{j+1}x, k_{l+1}x]+$
\begin{eqnarray}\label{anal2.5}(-1)^{(k+j)i}[h_k \cdot k_{l+1}x, g_{j+1}x,f_i]+
(-1)^{(l+j)k}[k_{l+1}x\cdot f_i, g_{j+1}x, h_k]=\end{eqnarray}
$(-1)^{(i+j)l+l+1}(\{f_ih_kg_{j+1},k_{l+1}\}-f_ih_k\{g_{j+1},k_{l+1}\})+\\
(-1)^{(k+j)i+j+1}(\{h_kk_{l+1},g_{j+1}\}f_i-\{h_kk_{l+1},g_{j+1}f_i\})+\\
(-1)^{(l+j)k+i+j+1}(\{k_{l+1}f_i,g_{j+1}\}h_k-\{k_{l+1}f_i,g_{j+1}h_k\}).$

Очевидно замечаем, что мы имеем аналог тождества (\ref{jorskob3}).

\medskip 

$0=
(-1)^{(i+j)l}[f_i \cdot h_{k+1}x, g_j, k_{l+1}x]+\\
(-1)^{(k+j)i}[h_{k+1}x \cdot k_{l+1}x, g_j,f_i]+
(-1)^{(l+j)k}[k_{l+1}x \cdot f_i, g_j, h_{k+1}x]=\\
(-1)^{(i+j)l+j+l+1}(\{f_ih_{k+1}g_j,k_{l+1}\}-\{f_ih_{k+1},g_jk_{l+1}\})+\\
(-1)^{(k+j)i+l+1}(\{h_{k+1},k_{l+1}\}g_jf_i-\{h_{k+1},k_{l+1}\}g_jf_i)+\\
(-1)^{(l+j)k+i+j+k+1}(\{k_{l+1}f_ig_j,h_{k+1}\}-\{k_{l+1}f_i,g_jh_{k+1}\}).$

Легко заметить, что полученное соотношение эквивалентно

$(-1)^{(i+j)l}(\{f_ih_kg_j,k_l\}-\{f_ih_k,g_jk_l\})=$
\begin{eqnarray}\label{anal2.6}
(-1)^{(l+j)k}(\{k_lf_ig_j,h_k\}-\{k_lf_i,g_jh_k\}).
\end{eqnarray}

\medskip 

$0=
(-1)^{(i+j)l}[f_i \cdot h_{k+1}x, g_{j+1}x, k_l]+\\
(-1)^{(k+j)i}[h_{k+1}x \cdot k_l, g_{j+1}x, f_i]+
(-1)^{(l+j)k}[k_l \cdot f_i, g_{j+1}x, h_{k+1}x]$

Заметим, что при замене $f_i \mapsto h_k, h_k \mapsto k_l, k_l \mapsto f_i$ мы имеем выражение (\ref{anal2.5}), 
которое эквивалентно (\ref{jorskob3}).

\medskip

$0=
(-1)^{(i+j)l}[f_{i+1}x\cdot h_k,g_j,k_{l+1}x]+\\
(-1)^{(k+j)i}[h_k\cdot k_{l+1}x, g_j, f_{i+1}x]+
(-1)^{(l+j)k}[k_{l+1}x\cdot f_ix,g_j,h_{k+1}]$

Легко заметить, что при замене $f_i \mapsto k_l, k_l \mapsto h_k, h_k \mapsto f_i$ мы получим выражение эквивалентное (\ref{anal2.6}).

\medskip

$0=
(-1)^{(i+j)l}[f_{i+1}x\cdot h_k, g_{j+1}x, k_l]+\\
(-1)^{(k+j)i}[h_k\cdot k_l, g_{j+1}x, f_{i+1}x]+
(-1)^{(l+j)k}[k_l\cdot f_{i+1}x, g_{j+1}x, h_k]$

Очевидно замечаем, что при замене $f_i \mapsto k_l, k_l \mapsto h_k, h_k \mapsto f_i$ мы имеем выражение (\ref{anal2.5}), 
которое эквивалентно (\ref{jorskob3}).
\medskip 

$0=
(-1)^{(i+j)l}[f_{i+1}x\cdot h_{k+1}x, g_j, k_l]+\\
(-1)^{(k+j)i}[h_{k+1}x \cdot k_l, g_j, f_{i+1}x]+
(-1)^{(l+j)k}[k_l \cdot f_{i+1}x, g_j, h_{k+1}x]$

Легко заметить, что при замене $f_i \mapsto h_k, h_k \mapsto k_l, k_l \mapsto f_i$ мы получим выражение эквивалентное (\ref{anal2.6}).

\medskip

Для доказательства теоремы осталось показать, что система тождеств (\ref{jorskob3},\ref{anal2.6}) 
эквивалентна системе тождеств (\ref{jorskob2},\ref{jorskob3}). 
Из (\ref{anal2.6}) и (\ref{jorskob3}) можем получить 
$$(-1)^{(i+j)l}\{f_ih_k,g_jk_l\}+(-1)^{(l+j)k}\{k_lf_ig_j,h_k\}+(-1)^{(k+j)i}\{h_kk_l,g_jf_i\}=$$
$$(-1)^{(i+j)l}f_ih_k\{g_j,k_l\}+(-1)^{(l+j)k}\{k_lf_i,g_j\}h_k+(-1)^{(k+j)i}\{h_kk_l,g_j\}f_i.$$
Откуда, путем замены $h_k \mapsto k_l, k_l \mapsto h_k$ получим 
$$(-1)^{(i+j)k}\{f_ik_l,g_jh_k\}+(-1)^{(k+j)l}\{h_kf_ig_j,k_l\}+(-1)^{(l+j)i}\{k_lh_k,g_jf_i\}=$$
$$(-1)^{(i+j)k}f_ik_l\{g_j,h_k\}+(-1)^{(k+j)l}\{h_kf_i,g_j\}k_l+(-1)^{(l+j)i}\{k_lh_k,g_j\}f_i.$$
Применим (\ref{jorskob3}) и получим
$$(-1)^{(i+j)k}f_ik_l\{g_j,h_k\}+(-1)^{(k+j)l}\{h_kf_i,g_j\}k_l+(-1)^{(l+j)i}\{k_lh_k,g_j\}f_i=$$
$$(-1)^{ik+il+lk}((-1)^{il+jl}f_ih_k\{g_j,k_l\}+(-1)^{(k+j)i}\{h_kk_l,g_j\}f_i+(-1)^{(l+j)k}\{k_lf_i,g_j\}h_k).$$
Откуда, путем замены $f_i \mapsto k_l, k_l \mapsto f_i$, получаем
$$(-1)^{(l+j)k}k_lf_i\{g_j,h_k\}+(-1)^{(k+j)l}\{h_kk_l,g_j\}f_i=$$
$$(-1)^{(k+j)i}h_kk_l\{g_j,f_i\}+(-1)^{(l+j)k}\{k_lf_i,g_j\}h_k.$$
Что, в свою очередь, является полным аналогом (\ref{jorskob2}). Таким образом, теорема доказана.

Отметим, что для неунитальной ассоциативно-суперкоммутативной супералгебры $\Gamma$ c дифференцированием $D$ 
супералгебра векторного типа $J(\Gamma,D)$ будет являться йордановой. Доказательство легко следует
из проверки тождеств (\ref{jorskob1}-\ref{jorskob3}).

\section{О $\delta$-супердифференцированиях обобщенного дубля Кантора с первичной четной частью}

Исследование $\delta$-дифференцирований вытекает из работ Н. С. Хопкинс \cite{Hop} и В. Т. Филиппова \cite{fi}, где рассматривались 
антидифференцирования (т.е. (-1)-дифференцирования) алгебр Ли. В дальнейшем, эти результаты получили обобщение в 
работах В. Т. Филиппова \cite{Fil,Fill}. Дальнейшее исследование $\delta$-дифференцирований связано с работами
В. Т. Филлипова \cite{ Filll}, И. Б. Кайгородова
\cite{kay,kay_lie, kay_lie2,kay_lie4}, И. Б. Кайгородова и В. Н. Желябина \cite{kay_lie3}, и П. Зусмановича \cite{Zus}. 
В результате, были описаны $\delta$-дифференцирования
первичных ассоциативных, альтернативных, лиевых и мальцевских нелиевых алгебр, полупростых йордановых алгебр, 
первичных лиевых супералгебр, полупростых конечномерных йордановых супералгебр и простых унитальных супералгебр йордановой скобки.
Более подробную информацию по описанию $\delta$-дифференцирований и $\delta$-супердифференцирований можно найти в 
обзоре \cite{kungur10}.

Однородный элемент $\psi$ суперпространства эндоморфизмов $A \rightarrow A$ называется
супердифференцированием, если
$$\psi(xy)=\psi(x)y+(-1)^{p(x)p(\psi)}x\psi(y).$$
Рассмотрим супералгебру Ли $A$ и зафиксируем элемент $x \in A_i$. Тогда $u_x: y \rightarrow [	x,y]$ 
является супердифференцированием супералгебры $A$ и его четность $p(u_{x})=i.$

Для фиксированного элемента $\delta$ из основного поля, под
$\delta$-супердиф\-фе\-рен\-ци\-ро\-ванием супералгебры $A=A_0\oplus A_1$ мы 
понимаем однородное линейное отображение $\phi : A \rightarrow A,$ такое что для однородных $x,y\in A$ выполнено
\begin{eqnarray}\label{der} \phi (xy)&=&\delta(\phi(x)y+(-1)^{p(x)p(\phi)}x\phi(y)).\end{eqnarray}

Под суперцентроидом $\Gamma_{s}(A)$ супералгебры $A$ мы будем понимать множество всех однородных линейных отображений $\chi: A \rightarrow A,$ для произвольных однородных элементов $a,b$ удовлетворяющих условию

$$\chi(ab)=\chi(a)b=(-1)^{p(a)p(\chi)}a\chi(b).$$

Центроид алгебры $A$ определяется по аналогии и обозначется $\Gamma(A).$

Заметим, что 1-супердифференцирование является обыкновенным супердифференцированием; 0-су\-пер\-диф\-фе\-ренци\-рованием
является произвольный эндоморфизм $\phi$ супералгебры $A$ такой, что $\phi(A^{2})=0$. 

Ненулевое $\delta$-супердифференцирование $\phi$ будем считать нетривиальным, если 
$\delta\neq 0,1$ и $\phi \notin \Gamma_s(A).$

В данной части мы рассмотрим $\delta$-супердифференцирования обобщенного дубля Кантора $J(\Gamma, \{, \})$, построенном исходя из 
$A=\Gamma$ --- первичной ассоциативной алгебры и скобки $\{,\}$. Напомним, что $A$ является первичной алгеброй, если из равенства $aAb=0$ для некоторых элементов $a,b$ алгебры $A$ следует, 
что либо $a=0$ либо $b=0.$

В дальнейшем, через $(\Delta_{\frac{1}{2}}(J(\Gamma, \{, \})))_i$ будем обозначать пространство $\frac{1}{2}$-супердифференцирований 
супералгебры $J(\Gamma, \{, \})),$ имеющих четность $i.$ 

\medskip

{\bf Теорема 3.1.} \emph{Супералгебра $J(\Gamma, \{, \})$ не имеет ненулевых $\delta$-супер\-дифференцирований, если $\delta\neq 0,\frac{1}{2},1.$
Если $\phi$ --- четное $\frac{1}{2}$-супер\-дифференци\-рование $J(\Gamma, \{, \})$, то 
$\{\phi|_A, \phi \in (\Delta_{\frac{1}{2}}(J(\Gamma, \{, \})))_0   \}=\Gamma(A) \cap \Delta_{\frac{1}{2}}(A, \{, \})$ 
и $\phi(ax)=\phi|_A(a)x.$ В частности, если $J(\Gamma, \{, \})$ ---
супералгебра векторного типа, то $\{\phi|_A, \phi \in (\Delta_{\frac{1}{2}}(J(\Gamma, \{, \})))_0   \}=\Gamma(A).$ }

\medskip 

{\bf Доказательство.} Пусть $\phi$ --- четное $\delta$-супердифференцирование и $\phi(cx)=\phi_c x.$ 
Рассмотрим ограничение $\phi|_A$ на подалгебру $A$. Таким образом, $\phi|_A$ есть $\delta$-дифференцирование 
первичной ассоциативной алгебры  $A$. Согласно \cite{Filll}, при $\delta=\frac{1}{2}$ имеем $\phi|_A \in \Gamma(A)$, а при 
$\delta\neq\frac{1}{2}$  имеем $\phi|_A=0.$ 

В дальнейшем, через $a,b,c,d$ мы обозначаем произвольные элементы алгебры $A$. Рассмотрим случай $\delta\neq\frac{1}{2},$ тогда 
$$\delta^{2}ab\phi(cx)=\delta a\phi(bcx)=\phi(abcx)=\delta ab\phi(cx),$$
что влечет $ab\phi(cx)=0,$ то есть $ab\phi_c=0.$ В силу первичности $A$, это влечет $\phi_c=0.$ Таким образом, 
мы имеем тривиальность $\phi.$

Если $\delta=\frac{1}{2},$ то
$$\frac{1}{2}\phi(a)bcx+\frac{1}{4}a\phi(b)cx+\frac{1}{4}ab\phi(cx)=\phi(abcx)=\frac{1}{2}\phi(ab)cx+\frac{1}{2}ab\phi(cx).$$
Откуда, пользуясь тем, что $\phi(a)bc=ab\phi(c)=a\phi(b)c=\phi(ab)c$, легко получаем $$ab(\phi_c-\phi(c))=0,$$ что, в силу первичности $A$, дает $\phi(c)=\phi_c.$
Отметим, что 
$$\phi \{a, b\}=\phi(ax \cdot bx)=\frac{1}{2}(\phi(ax)\cdot bx+ax \cdot \phi(bx))=\frac{1}{2}\{\phi(a),b\}+\frac{1}{2}\{a,\phi(b)\},$$
то есть $\phi|_A \in \Delta_{\frac{1}{2}}(A,\{,\}).$ Заметим, что любой элемент
$\psi|_A \in \Gamma(A) \cap \Delta_{\frac{1}{2}}(A, \{, \})$ продолжается до четного $\frac{1}{2}$-супердифференцирования $\psi$
по правилу $\psi(ax)=\psi|_A(a)x.$

Пусть $\phi$ --- нечетное $\frac{1}{2}$-супердифференцирование супералгебры $J(\Gamma, \{, \}),$ тогда
$$\delta^{2}\phi(a)bc+\delta^{2}a\phi(b)c+\delta ab\phi(c)=\phi(abc)=\delta \phi(a)bc+\delta^2a\phi(b)c+\delta^2ab\phi(c).$$

Откуда легко следует $\phi(a)bc=ab\phi(c),$ пользуясь чем, имеем
$$\phi(ab)dc=abd\phi(c)=a\phi(b)dc.$$
Откуда, в силу первичности $A,$ вытекает $\phi(ab)=\phi(a)b=a\phi(b).$ Воспользовавшись (\ref{der}), имеем
$$\phi(ab)=\delta(\phi(a)b+a\phi(b))=2\delta\phi(ab).$$
Таким образом, возможны два случая $\delta=\frac{1}{2}$ или  $\delta\neq \frac{1}{2}$ и $\phi(ab)=0.$ 

Из второго случая, легко вытекает $0=\phi(ab)c=\phi(a)bc.$ Откуда в силу первичности $A,$ получаем $\phi(A)=0.$
Осталось заметить, что 
$$\delta ab\phi(cx)+\delta\phi(ab)\cdot cx=\phi(abcx)=\delta\phi(a)(bcx)+\delta^2a(\phi(b)\cdot cx)+\delta^2ab\phi(cx),$$
то есть $(\delta^2-\delta)ab\phi(cx)=0$. Таким образом, используя первичность $A$, имеем $\phi(cx)=0.$ 
Основное утверждение теоремы доказано.

Отметим, что если $J(\Gamma, \{, \})$ --- супералгебра векторного типа и $\phi|_A \in \Gamma(A),$ то 
$$D(a\phi|_A(b)-\phi|_A(a)b)=0,$$
$$D(a)\phi|_A(b)+aD(\phi|_A(b))-(D(\phi|_A(a))b-\phi|_AD(b))=0,$$
$$2\phi|_A(D(a)b-aD(b))=D(\phi|_A(a))b-\phi|_A(a)D(b)+D(a)\phi|_A(b)-aD(\phi|_A(b)),$$
$$\phi|_A\{a,b\}=\frac{1}{2}\{\phi|_A(a),b\}+\frac{1}{2}\{a,\phi|_A(b)\},$$
то есть $\phi|_A \in \Delta_{\frac{1}{2}}(A,\{,\}).$
Теорема доказана.

\medskip

Отметим, что в случае когда $\Gamma$ --- является ассоциативно-коммутативной первичной алгеброй и скобка $\{,\}$ является 
йордановой скобкой, то по теореме 2.1 мы получаем отсутствие ненулевых $\delta$-супердифференцирований 
при $\delta \neq 0,\frac{1}{2},1$
для супералгебр йордановой скобки с первичной коммутативно-ассоциативной четной частью.

\end{document}